\documentclass{elsart}
\usepackage{amscd,amssymb,amsmath,amsfonts,latexsym}
\usepackage[arrow,matrix,graph,frame,poly,arc,tips]{xy}
\journal{Journal of Algebra}

\newtheorem{theorem}[subsection]{Theorem}
\newtheorem{lemma}[subsection]{Lemma}
\newtheorem{proposition}[subsection]{Proposition}
\newtheorem{corollary}[subsection]{Corollary}
\newtheorem{example}[subsection]{Example}

\newtheorem{mthm}{Theorem} 
\numberwithin{equation}{section}

\newcommand{\abs}[1]{\left|#1\right|}

\newcommand{\DS}{\displaystyle }
\newcommand{\surj}{\twoheadrightarrow}

\newcommand{\nor}{\triangleleft}

\newcommand{\C}{\mathbb{C}}

\newcommand{\Z}{\mathbb{Z}}

\renewcommand{\P}{P}

\newcommand{\g}{{\gamma}}
\newcommand{\G}{\Gamma}
\newcommand{\s}{\sigma}
\renewcommand{\a}{{\alpha }}
\renewcommand{\b}{{\beta}}
\renewcommand{\l}{\lambda}
\renewcommand{\t}{\tau}

\newcommand{\rank}{\operatorname{rank}}
\newcommand{\corank}{\operatorname{corank}}

\newcommand{\id}{\operatorname{id}}
\newcommand{\Aut}{\operatorname{Aut}}
\newcommand{\Hom}{\operatorname{Hom}}
\newcommand{\End}{\operatorname{End}}
\newcommand{\Epi}{\operatorname{Epi}}

\newcommand{\GL}{\operatorname{GL}}
\newcommand{\SL}{\operatorname{SL}}

\newcommand{\ab}{\operatorname{ab}}

\def\Im{\operatorname{Im}}
\def\im{\operatorname{im}}

\newcommand{\BS}{\operatorname{BS}}

\newcommand{\PSp}{\operatorname{PSp}}

\newcommand{\md}[1]{\ (\bmod\:{#1})}

\begin{document}
\begin{frontmatter}

\title{Counting homomorphisms onto finite \\ solvable groups}

\author[imar,nu]{Daniel~Matei}
\ead{dmatei@styx.math.neu.edu}
\author[nu]{Alexander~I.~Suciu\thanksref{nsf}}
\ead{a.suciu@neu.edu}
\thanks[nsf]{Supported by NSF grant DMS-0311142}

\address[imar]{Institute of Mathematics of the Academy,
P.O. Box 1-764, RO-014700 Bucharest, Romania}
\address[nu]{Department of Mathematics, Northeastern University,
Boston, MA 02115}

%\ead[url]{http://www.math.neu.edu/matei}
%\ead[url]{http://www.math.neu.edu/\~{}suciu}


%\date{April 20, 2004}
%\date{July 21, 2004}
%\date{October 30, 2004}
%\date{December 14, 2004}
%\date{January 8, 2005}

\begin{abstract}
We present a method for computing the number of epimorphisms 
from a finitely presented group $G$ to a finite solvable group $\G$, 
which generalizes a formula of Gasch\"{u}tz.  
Key to this approach are the degree $1$  and $2$ cohomology 
groups of $G$, with certain twisted coefficients.  As an application, 
we count low-index subgroups  of $G$. We also investigate the finite
solvable quotients of the Baumslag-Solitar groups, the Baumslag 
parafree groups, and the Artin  braid groups.
\end{abstract}

\end{frontmatter}

\section{Introduction}
\label{sec:intro}

\subsection{Finite quotients}
\label{subsec: finite quotients}
One way to understand an infinite, finitely generated group is
to identify its finite quotients, and count all the epimorphisms 
to one of these finite groups.  A wide spectrum of possibilities can 
occur. For example, residually finite groups have plenty 
of finite quotients, whereas infinite simple groups have none.  
Free groups and surface groups have an 
abundance of finite solvable quotients, whereas groups 
with perfect derived subgroup have no solvable quotients 
except abelian ones.

If $G$ is a finitely generated group, and $\G$ a finite group, 
let $\abs{\Hom(G,\G)}$ be  the number 
of homomorphisms from $G$ to $\G$, and let 
$\delta_{\Gamma}(G)=\abs{\Epi(G,\G)}/\abs{\Aut \G}$ be 
the number of epimorphisms from $G$ to $\G$, up to 
automorphisms of $\G$. In the case when $G=F_n$ is the free 
group of rank $n$, Philip Hall   \cite{HaP} gave a procedure 
to determine the Eulerian function $\abs{\Epi(F_n,\G)}$, 
based on M\"obius inversion in the subgroup lattice of $\G$. 
An explicit formula for computing $\abs{\Epi(F_n,\G)}$ in 
the case when $\G$ is solvable was given by Gasch\"{u}tz \cite{Ga}.  

In this paper, we generalize Gasch\"{u}tz's formula, from the 
free group $F_n$, to an arbitrary finitely presented group $G$.  
As a byproduct, we derive an expression for the order 
of the automorphism group of a finite solvable group $\G$. 
Putting things together gives a method for computing the 
solvable Hall invariants $\delta_{\Gamma}(G)$ 
in terms of homological data.  
This extends previous results from \cite{MS}, 
which only dealt with certain  metabelian groups $\G$. 

\subsection{Finite-index subgroups}
\label{subsec:finite index}
Another way to understand a finitely generated, residually finite 
group is through its finite-index subgroups. Let $a_k(G)$ and 
$a_k^{\nor}(G)$ be the number of index $k$ subgroups 
(respectively, normal subgroups) of $G$. The growth of 
these sequences---also known as the subgroup growth of 
$G$---has been a subject of intensive study in the recent past, 
see \cite{LS}.   Much is known in the case when $G$ is nilpotent; 
explicit formulas for $a_k(G)$ and $a_k^{\nor}(G)$ 
are available in a few other cases, such as free products 
of cyclic groups and surface groups. 

In \cite{Ha}, Marshall Hall 
showed how to express the numbers $a_k(G)$ in terms of the Hall 
invariants $\delta_{\G}(G)$, where $\G$ ranges through the 
isomorphism classes of subgroups of the symmetric group $S_k$. 
In \cite{MS}, we used this fact to arrive at a homological formula 
for $a_3(G)$.  Here, we give a similar (but more involved) formula 
for $a_4(G)$. Combining our previous results with the present 
techniques, we also give formulas for $a_k^{\nor}(G)$,  for $k\le 15$. 

\subsection{Solvable quotients}
\label{subsec:solv quot}
The derived series of a group $G$ is defined inductively by 
$G^{(0)}=G$ and $G^{(k)}=[G^{(k-1)},G^{(k-1)}]$.
A group $G$ is solvable if its derived series terminates.  
The derived length of $G$ is the minimal $k$ for which $G^{(k)}=1$; 
abelian groups have derived length $1$, while metabelian groups 
have length $2$.  

At the other extreme, a perfect group $G$ equals its own derived 
subgroup, so its derived series stabilizes at $G'=G^{(1)}$.  
Clearly, a solvable group has no perfect subgroups. Hence, 
if $G^{(k)}$ is perfect, then $G$ has 
no solvable quotients of derived length greater than $k$.

Now suppose $\G$ is a finite solvable group, of derived length $l$.  
Since the derived series of $G$ consists of characteristic subgroups, 
the number of epimorphisms from $G$ to $\G$ depends only 
on the solvable quotient $G/G^{(l)}$; 
in fact, $\abs{\Epi(G,\Gamma)}=\abs{\Epi(G/G^{(l)},\Gamma)}$.

\subsection{Lifting homomorphisms}
\label{subsec: lifting}
As is well-known, a group  is solvable if and only if it can 
be expressed as an iterated extension of abelian 
groups. In order to count homomorphisms from a finitely generated 
group $G$ to a finite solvable group $\G$, we use an inductive 
procedure, sketched below.   

Suppose we have an extension $1\to A\to\G\to B\to 1$, with $A$ 
abelian and $B$ solvable.  Such an extension is 
determined by a monodromy homomorphism $\s\colon B\to \Aut(A)$,  
and a (twisted) cohomology class $[\chi]\in H_{\s}^2(B, A)$. Let 
 $\rho\colon G\to B$ be a homomorphism.  Then $\rho$ 
has a lift $\tilde\rho\colon G\to \G$ if and only if 
$\rho^*[\chi]=0$ in $H_{\s\rho}^2(G, A)$. Furthermore, 
the lifts of $\rho$ are in one-to-one correspondence with 
$1$-cochains that cobound the $2$-cocycle $-\rho^*\chi$. 
The set  of such $1$-cochains, $Z_{\s\rho,\chi}^1(G, A)$, is 
either empty, or is in bijection with $Z_{\s\rho}^1(G, A)$; 
define $\epsilon_{\chi}(\rho)$ to be  $0$ or $1$, accordingly.  
We find:
\[
\abs{\Hom(G,\G)}= \sum_{\rho\in \Hom(G,B)} 
\epsilon_{\chi}(\rho) \cdot \abs{Z_{\s\rho}^1(G, A)}.
\]

\subsection{Systems of equations}
\label{subsec: systems}
If $G$ admits a finite presentation, say 
$\mathcal{P}=\langle x_1,\dots ,x_n\mid r_1,\dots ,r_m\rangle$, 
we can translate the lifting condition into a system 
$S=S(\mathcal{P},\G, \rho)$ of $m$ equations in $n$ unknowns 
over the abelian group $A$.  The equations in $S$, given in 
(\ref{eq:foxnh}) below, are similar in nature to inhomogeneous 
linear equations. The homogeneous part is written in terms 
of the Fox derivatives $\partial r_i/\partial x_j$, 
twisted  by $\s\rho$, whereas the non-homogeneous 
part involves only $\rho^{*}\chi$ and the presentation $\mathcal{P}$.

As shown in  Theorem \ref{thm:fpcount}, 
the number of solutions of the system $S$  coincides 
with the number of lifts of $\rho$. The precise determination
of these solutions gives a way to explicitly construct those lifts.

In general, the system of equations $S$ cannot be reduced
to a linear system. But, in certain cases, this is possible. One instance 
(exploited in \cite{Fx2} and \cite{MS}) is when $B$ is abelian, 
$A$ is the additive group of a finite commutative ring $R$, 
and $\s$ is induced by multiplication 
in $R$. Another instance is when $A$ is homocyclic, say
$A=\Z_{q^r}^{\oplus s}$, and the monodromy is of 
the form $\s\colon B\to \GL(s,\Z_{q^r})$.
In particular, if $A$ is an elementary abelian $q$-group $E$, 
the cohomology group $H^1(G, E)$ can be viewed as a vector 
space over the prime field $\Z_q$, and its dimension can be 
computed in terms of the rank of the Jacobian matrix associated 
to a presentation of $G$, twisted by the homomorphism $\rho$. 

\subsection{A generalization of Gasch\"{u}tz' formula}
\label{subsec: gengas}
Next, we restrict our attention to surjective homomorphisms 
$G\surj \G$. From the above discussion, we know that 
$\abs{\Hom(G,\G)}$ can be obtained by summing 
$\epsilon_{\chi} (\rho) q^{\dim_{\Z_q}Z^{1}_{\s\rho}(G,E)}$ 
over all $\rho\in \Hom(G,B)$.   In order to compute 
$\abs{\Epi(G,\G)}$, we have to subtract all the 
homomorphisms $G\to\G$ which are not surjective. 
Due to the minimality of $E$, those are precisely the 
homomorphisms whose image is a complement of $E$ in $\G$. 

Every solvable group $\G$ admits a normal series whose successive 
factors are elementary abelian. Our method for calculating 
$\abs{\Epi(G,\Gamma)}$ is to use such a chief series
to construct all the epimorphisms by repeated liftings through 
the chief series.  The key step is provided by the following result. 

\begin{mthm}
\label{thm:gg intro}
Let $G$ be a finitely presented group, and let $\G$ be a finite, 
solvable group.  Let $E$ be an elementary abelian $q$-group 
which is also a chief factor of $\G$, and let $B=\G/E$. If $\s$ is 
the monodromy and $\chi$ is the $2$-cocycle defining
the extension $1\to E\to\G\to B\to 1$, then:
\[
\abs{\Epi(G,\G)}=\abs{E}^{\zeta}  \sum_{\rho\in\Epi(G,B)}
\left( \epsilon_{\chi} (\rho) q^{\dim_{\Z_q}H^{1}_{\s\rho}(G,E)} - 
c_{\chi} q^{\kappa(\a-1)} \right), 
\]
where $\zeta=0$ or $1$ according as $B$ acts 
trivially on $E$ or not, 
$\epsilon_{\chi}(\rho)=1$ or $0$ according as the equation  
$\delta^1 f=-\rho^*\chi$ has a solution  or not,
$c_{\chi}=1$ or $0$ according as 
$[\chi]\in H^2_{\s}(B,E)$ vanishes or not, 
$\a$ is the number of complemented chief factors of $\G$ 
isomorphic to $E$ as $\G$-modules (under the conjugation action), 
and $q^{\kappa}=\abs{\End_{\G}(E)}$.  
\end{mthm}

Now suppose $\G=\G_0>\G_{1}> \cdots> \G_{\nu}> \G_{\nu+1}=1$ 
is a chief series, with factors $E_i=\G_i/\G_{i+1}=\Z_{q_i}^{\oplus s_i}$ 
and quotients $B_i=\G/\G_{i}$.  We then find:
\[
\abs{\Epi(G,\G)}=
\sum_{\rho_{1}\in\Epi_{\rho_{0}}(G,B_{1})} \!\cdots\!
\sum_{\rho_{\nu}\in\Epi_{\rho_{\nu-1}}(G,B_{\nu})}
\abs{E_{\nu}}^{\zeta_{\nu}} \! \left( \epsilon_{\chi_{\nu}}(\rho_{\nu})  
q_{\nu}^{\beta_{\nu}} - 
c_{\chi_{\nu}} q_{\nu}^{\kappa_{\nu}(\a_{\nu}-1)} \right) \! ,
\]
where  $B_{\nu+1}=E_{\nu}\times _{\s_\nu, \chi_\nu} B_{\nu}$, 
$\beta_{\nu}=\dim_{\Z_{q_{\nu}}} 
H^{1}_{\s_{\nu}\rho_{\nu}} (G, E_{\nu})$, 
let $\a_\nu$ be the number of chief factors 
of $B_{\nu+1}$ isomorphic to $E_\nu$ as $B_{\nu+1}$-modules, 
and $\Epi_{\rho_{i}}(G,B_{i+1})$ is 
the set of epimorphisms lifting $\rho_{i}\colon G\surj B_{i}$.
In the case when $G=F_n$, this recovers Gasch\"{u}tz' formula. 

\subsection{Examples}
\label{subsec: examples}
To illustrate our recursive process for calculating $\abs{\Epi(G,\Gamma)}$, 
we discuss various classes of source and target groups.  

When it comes to the target group $\G$, we analyze in detail two 
series of finite metabelian groups: the dihedral and the binary dihedral groups. 
We also consider a class of derived length $3$ solvable groups, 
of the form $\G=\Z_q^2\rtimes D_{2p}$, which includes the symmetric group 
$S_4=\Z_2^{\oplus 2}\rtimes D_6$.

In calculating $\abs{\Epi(G,\Gamma)}$, one can use the lattice of 
subgroups of $\G$ instead of its chief series extensions.
We briefly illustrate this approach for the sake of comparison, 
in the case when $\G$ is a dihedral group.

When it comes to the source group $G$, we start of course 
with the free groups $F_n$.  Another family of examples are the 
 orientable and non-orientable surface groups, 
$\Pi_g$ and $\Pi_g^*$.  The other examples we consider 
(the one-relator Baumslag-Solitar and Baumslag groups, a certain 
link group, and the Artin braid groups) are discussed separately below. 

\subsection{Baumslag-Solitar groups}
\label{subsec: bsg}
A famous family of one-relator groups 
was introduced by Baumslag and Solitar in \cite{BS}.  
For each pair of integers $(m,n)$ with  $0<m\le \abs{n}$, let 
$\BS(m,n)=\langle x,y \mid x y^{m} x^{-1} y^{-n} \rangle$. 
Much is known about these groups:  
$\BS(m,n)$ is solvable if and only if   $m=1$, 
in which case $\BS(1,n)=\Z[1/n]\rtimes \Z$; it is residually finite 
if and only if  $m=\abs{n}$ or $m=1$, in which  
case it is also Hopfian; and it is Hopfian
if and only if  $m$ and $n$ have the same prime divisors or 
$m=1$.  For example, $\BS(2,3)$ is 
non-Hopfian, while $\BS(2,4)$  is 
Hopfian but non-residually finite. 

The groups $\BS(m,n)$ have been classified by Moldavanski 
\cite{Mo}.  They are in bijection with the set of unordered pairs 
$(m, n)$ with $0<m\le \abs{n}$.  In the case when $m=1$, 
the set of finite quotients of $\BS(1,n)$ is a complete group 
invariant, see \cite{MoS}, and it consists of all quotients
of metacyclic groups of type $\Z_s\rtimes_{\s}\Z_r$, where 
$\s$ is multiplication by $n$, and $n^r\equiv 1 \md{s}$. 
The subgroup growth of the Baumslag-Solitar groups 
$\BS(m,n)$, with $m$, $n$ coprime was 
determined by E.~Gelman, see \cite[p.~284]{LS}; 
a presentation for $\Aut(\BS(m,n))$ was given in \cite{GHMR}.

We compute here the number of epimorphisms from the 
Baumslag-Solitar groups to $D_8$ and $Q_8$. 
This allows us to divide the groups $\BS(m,n)$ into four 
(respectively, two) non-isomorphic classes.

\subsection{Baumslag's parafree groups}
\label{subsec: bpg}
In \cite{Ba}, Baumslag introduced the following notion:
A group $G$ is called parafree if it is  residually nilpotent and 
has the same nilpotent quotients as a free group $F$.  The simplest 
non-free yet parafree groups are the one-relator groups 
$\P(m,n)=\langle x,y,z \mid x z^{m} x z^{-m} x^{-1} z^{n} y z^{-n} y^{-1} \rangle$. 
As shown in  \cite{Ba}, each group in this family (indexed by pairs of 
integers $m$ and $n$) has the same nilpotent quotients 
and the same first two solvable quotients as the free group $F_2$. 

In \cite{LL},  R. Lewis and S. Liriano showed that there are several 
distinct isomorphism types among the groups $\P(m,n)$.  
By counting homomorphisms to $\SL(2,\Z_4)$, they verified that 
the third solvable quotients of $\P(m,n)$ differ from those of $F_2$, 
for certain pairs of integers, e.g., $(m,n)=(1,3)$ and $(3,5)$. 
By computing the number of epimorphisms of $\P(m,n)$ onto a 
smaller solvable quotient group of derived length $3$, namely $S_4$, 
we can recover (and sharpen) the result of Lewis and Liriano.  
We find:  If $m$ odd and $m-n \equiv 2\md{4}$, then $\P(m,n)$ 
is not isomorphic to $F_2$. Moreover, in that case, $\P(m,n)$ is not
isomorphic to $\P(m',n')$ if $m'$ even or $m'-n' \not\equiv 2 \md{4}$.

\subsection{A link group}
\label{subsec: link}
Counting finite solvable quotients of a group $G$ can also help 
decide whether a normal subgroup $H$ is perfect.
Indeed, if $H\nor G$ is perfect, and $\G$ is finite and solvable, 
then $\delta_{\G}(G)=\delta_{\G}(G/H)$. In other words, if 
$\delta_{\G}(G)>\delta_{\G}(G/H)$ for some finite solvable 
group $\G$, then $H$ is not perfect. 

As an example of how this works, 
we take $G$ to be the group of a certain $2$-component link 
considered by Hillman in \cite{Hil}. Then $G_{\omega}$, the 
intersection of the lower central series of $G$, is non-trivial, i.e., 
$G$ is not residually nilpotent. Hence, $G$ is not parafree.  
On the other hand, $G/G_{\omega} \cong \P(-1,1)$, and thus 
$G$ has the same nilpotent quotients as $F_2$; moreover, 
$G/G''\cong F_2/F_2''$. Using our formula, we find that 
$\delta_{S_4}(G)>\delta_{S_4}(\P(-1, 1))$.  
This shows that $G_{\omega}$ is not perfect, thereby 
answering a question of Hillman, see \cite[p.~74]{Hil}.  

\subsection{Braid groups}
\label{subsec: artin}
The braid groups $B_n$  have been intensively studied ever 
since Artin introduced them in the mid 1920's.  
Although the braid groups are  residually finite, 
few of their finite quotients are known. 
General series of non-abelian quotients for $B_n$ ($n\ge 3$) 
are the symmetric groups $S_n$, and 
the projective symplectic groups $\PSp(n-2, 3)$ with $n$ even, 
or $\PSp(n-1, 3)$ with $n$ odd. 

Once we restrict to solvable quotients, the situation becomes 
more manageable.  In Section \ref{sec:braids}  
we use our methods to compute the number of epimorphisms 
from $B_n$ to certain finite solvable groups.  For example, we 
show that $\delta_{S_3}(B_3)=1$ and $\delta_{S_4}(B_4)=3$ 
(this recovers a particular case of a much more general result 
of Artin \cite{Ar}, see also \cite{CP}). Using results of V. Lin, 
we also compute the number of index $k$ subgroups of $B_n$, 
when $k\le n$ or $k=2n$, and $n$ is sufficiently large. 

We conclude with some conjectures on the possible values for 
$\delta_{\G}(B_n)$, for solvable $\G$, and on the behavior of the 
sequence $a_k(B_n)$, for $n\gg 0$.    

\section{Extensions and group cohomology}
\label{sec:homology}

In this section, we review some basic material on group cohomology. 
We outline a computation method based on Fox calculus, and 
explain how low-degree cohomology is connected with extensions
with abelian kernel.  We use \cite{Br82} and \cite{HS} as general reference. 

\subsection{Group cohomology and Fox calculus} 
\label{subsec:fox}

Let $G$ be a group, and $A$ a $G$-module, with action specified 
by a homomorphism $\alpha\colon G \to \Aut(A)$. Let 
$C^r=\operatorname{Maps}(G^{\times r},A)$ be the group of $r$-cochains,  
and define coboundary maps $\delta^r\colon C^r\to C^{r+1}$ 
by $\delta^{r}(f)(x_0,\dots ,x_r)=x_0 f(x_1,\dots ,x_r)- 
\sum_{i=0}^{r-1} (-1)^i f(x_0,\dots , x_ix_{i+1},\dots ,x_r)+
(-1)^{r-1}f(x_0,\dots, x_{r-1})$.  The cohomology groups of $G$ 
with coefficients in $A$ are defined as 
\begin{equation}
\label{eq:coho1}
H^r_{\a}(G,A)=Z^r_{\a}(G,A)/B^r_{\a}(G,A), 
\end{equation}
where $Z^r_{\a}(G,A)=\ker (\delta^r)$ are the cocycles 
and $B^r_{\a}(G,A)=\im (\delta^{r-1})$ are the coboundaries.  

Now suppose $G$ admits a finite presentation, $G=F_n/R$, 
where $F_n$ is the free group on $x_1,\dots, x_n$ and 
$R$ is the normal subgroup generated by $r_1,\dots, r_m$. 
The Fox derivatives 
$\frac{\partial}{\partial x_j}\colon \Z F_n\to\Z F_n$ are 
the $\Z$-linear maps defined by 
$\frac{\partial x_i}{\partial x_j}=\delta_{ij}$ and 
$\frac{\partial (uv)}{\partial x_j}=\frac{\partial u}{\partial x_j}
\epsilon(v) +u\frac{\partial v}{\partial x_j}$, 
where $\epsilon\colon \Z F_n\to \Z$ is the augmentation map. 
The beginning of a free resolution of $\Z$ by $\Z G$-modules 
is then 
\begin{equation}
\label{eq:foxres}
\xymatrix{
\Z G^{m} \ar[r]^{J_G} & \Z G^{n}  \ar[r]^{d_1}
& \Z G  \ar[r]^{\epsilon} &\Z \ar[r] & 0
},
\end{equation}
where $d_1=\phi (x_1-1 \,\dots\,  x_n-1)^{\top}$  and $J_G=
\phi \left(\partial r_i / \partial x_j\right)$ is the Fox Jacobian matrix. 
Applying $\Hom_{\Z G}(-,A)$ yields the cochain complex 
\begin{equation}
\label{eq:chain1}
\xymatrix{A \ar[r]^{d_1^{\a}}  & A^{n} 
\ar[r]^{J_G^{\a}} & A^{m}},
\end{equation}
whose homology is $H^1_{\a}(G,A)$. 

\subsection{Extensions with abelian kernel}
\label{subsec:abelker}
An extension $\G$ of a group $B$ by an abelian group $A$ 
(written additively) is a short exact sequence 
\begin{equation}
\label{eq:ext}
\xymatrix{1\ar[r] & A \ar[r]^{i} & \G\ar[r]^{\pi} & B\ar[r] & 1},
\end{equation}
determined by:
\begin{itemize}
\item the monodromy homomorphism $\s\colon B\to \Aut(A)$, 
defined by $i(\s_{b}(a))=s(b)\cdot i(a) \cdot s(b)^{-1}$, 
where $s\colon B\to\G$ is any set section of $\pi$;
\item the cohomology class $[\chi]\in H^2_{\s}(B,A)$ of a 
normalized $2$-cocycle $\chi \colon B\times B\to A$, defined by 
$i(\chi(b,b'))=s(b)s(b')s(bb')^{-1}$. 
\end{itemize}
An element of $\G$ can we written as a pair $(a, b)$ 
with $a\in A$ and $b\in B$, while multiplication in $\G$ is given by 
$(a_1, b_1)\cdot (a_2, b_2)=(a_1+\s_{b_1}(a_2)+\chi(b_1, b_2), b_1 b_2)$. 
Note that $(a, b)^{-1}=(-\s_{b^{-1}}(a)-\chi(b^{-1}, b), b^{-1})$.

We denote a group $\Gamma$ as in (\ref{eq:ext}) by 
$\Gamma=A \times_{\sigma, \chi} B$. In the case of a split extension 
($[\chi]= 0$), we simply write $\Gamma=A\rtimes_{\s} B$; in the case 
of a central extension ($\s_b = \id$, for all $b\in B$), 
we write $\Gamma=A\times_{\chi} B$.

Now assume $\G=A \times_{\sigma, \chi} B$ is finite, and let 
$c(\G)$ be the number of complements of $A$ in $\G$.  If the 
extension splits (which we write as $c_{\chi}=1$), then 
$c(\G)=\abs{Z^1_{\s}(B,A)}$; if the extension does not split 
(which we write as $c_{\chi}=0$), then $c(\G)=0$.  Thus: 
\begin{equation}
\label{eq:complements}
c(\G)=c_{\chi} \abs{Z^1_{\s}(B,A)}. 
\end{equation}

\section{Lifting homomorphisms}
\label{sec:lift}

In this section, we present a method for counting 
the homomorphisms from a finitely presented group $G$ 
to a finite group $\G$, given as an extension 
with abelian kernel.

\subsection{Homomorphisms into extensions}
\label{subsec:counthom}
Let $G$ be a finitely generated group, and $\rho\colon G\to B$ 
a homomorphism to a group $B$.   Let 
$\Gamma=A\times_{\s,\chi} B$ be an extension of $B$ by an 
abelian  group $A$, with monodromy $\s\colon B\to \Aut(A)$ 
and $2$-cocycle $\chi\colon B \times B \to A$. Then $A$ 
becomes a $G$-module, with action specified by 
$\s\rho \colon G\to \Aut(A)$.  Let 
$Z^1_{\s\rho} (G,A)$  and $H^1_{\s\rho} (G,A)$ be the 
corresponding $1$-cocycle and $1$-cohomology groups. 

Any set map $\l\colon G\to\G$ lifting $\rho\colon G\to B$ 
can be written as a pair of maps $\l=(f, \rho)$, with 
$f\colon G\to A$.  Then $\l$ is a homomorphism if and 
only if $f$ satisfies the following: 
\begin{equation}
\label{eq:homcond}
f(gh)=f(g)+\s_{\rho(g)}(f(h))+\chi(\rho(g),\rho(h)), 
\quad \hbox{for all $g,h\in G$}, 
\end{equation}
that is, $f$ is a $1$-cochain that cobounds the $2$-cocycle 
$- \rho^*\chi \colon G\times G \to A$.  Denote the set of all 
such $1$-cochains by 
\begin{equation}
\label{eq:z1chi}
Z^1_{\s\rho,\chi } (G,A) = \{ f\colon G\to A \mid \delta^1 f = 
- \rho^* \chi \}. 
\end{equation}
This set is either empty, or else 
$Z^1_{\s\rho,\chi } (G,A)=Z^1_{\s\rho } (G,A) + f_0$, for some 
$1$-cochain $f_0$. Set $\epsilon_{\chi}(\rho)$ to be $0$ or $1$, 
accordingly.  We then have the following. 

\begin{proposition}
\label{prop:homcount}
The number of homomorphisms from $G$ to 
$\G=A\times_{\s,\chi} B$ is given by
\begin{equation}
\label{eq:homchi}
\abs{\Hom(G,A\times_{\s,\chi} B)}= 
\sum_{\rho\in\Hom(G,B)} \epsilon_{\chi}(\rho) \cdot 
\abs{Z^1_{\s\rho} (G,A)}. 
\end{equation}
\end{proposition}

If the extension is split, then clearly $\epsilon_0(\rho)=1$ 
for all $\rho\colon G\to B$, and thus (\ref{eq:homchi}) reduces to 
$\abs{\Hom(G,A\rtimes_{\sigma} B)}= \sum_{\rho} \abs{Z^1_{\s\rho} (G,A)}$. 
If the extension is central, then $Z^1_{\s\rho} (G,A)=\Hom(G,A)$ 
for all $\rho$, and thus $\abs{\Hom(G,A\times_{\chi} B)}= 
\sum_{\rho} \epsilon_{\chi}(\rho) \abs{\Hom (G,A)}$. 
If the extension is both split and central (i.e., a direct product), then 
 (\ref{eq:homchi}) reduces to the well-known formula 
$\abs{\Hom(G,A\times B)}= \abs{\Hom (G,A)} \cdot \abs{\Hom (G,B)}$. 

\subsection{Equations for lifts}
\label{subsec:eqlift}
We now give a practical algorithm for computing the 
quantities involved in Formula (\ref{eq:homchi}), 
in the case when $G$ is a finitely presented group.  
Given a homomorphism 
$\rho\colon G \to B$, we want to decide whether there is 
a $1$-cochain $f\colon G\to A$ cobounding $-\rho^*\chi$ 
(i.e., whether $\epsilon_{\chi}(\rho)\ne 0$), 
and, if so, count how many such cochains there are (i.e., determine 
$\abs{Z^1_{\s\rho} (G,A)}$).

Let $\mathcal{P}=\langle x_1,\dots ,x_n \mid r_1, \dots ,r_m \rangle$ 
be a finite presentation for $G$, and $\phi\colon F_n \to G$ the 
presenting homomorphism.  Write $r_k=u_{k,1} \cdots u_{k,l_k}$, 
with each $u_{k,j}$ equal to some $x_{i}^{e_{k,j}}$, where 
$e_{k,j}=\pm 1$.  

\begin{theorem}
\label{thm:fpcount}
Let $\rho\colon G \to B$ be a homomorphism.  
Then $\epsilon_{\chi}(\rho)=1$ or $0$ according to whether 
the following system of equations has a solution 
$(a_1,\dots ,a_n)$ with  $a_i\in A$:
\begin{equation}
\label{eq:foxnh}
\begin{split}
\sum_{i=1}^{n} \s\bar\rho 
\bigg(\frac{\partial r_k}{\partial x_i}\bigg) & (a_i)  +
\sum_{j=1}^{l_k}  \frac{e_{k,j}-1}{2}\chi 
\big( \bar\rho (u_{k,j} ), \bar\rho(u_{k,j}^{e_{k,j}}) \big) 
\\
& + \sum_{j=1}^{l_k-1}  \chi \big(
\bar\rho  (u_{k,1} \cdots  u_{k,j} ),
\bar\rho (u_{k,j+1} ) \big) =0, \quad 1\le k\le m,
\end{split}
\end{equation}
where $\bar\rho=\rho \phi$.  Moreover, if $\epsilon_{\chi}(\rho)=1$, then 
$\abs{Z^1_{\s\rho} (G,A)}$ equals the number of solutions 
of the (homogeneous) system
\begin{equation}
\label{eq:foxhom}
\sum_{i=1}^{n} \s\bar\rho\big(\tfrac{\partial r_k}{\partial x_i}\big)(a_i) =0, 
 \qquad 1\le k\le m.
\end{equation}
\end{theorem}

\begin{pf}
Let $f\colon G\to A$ be a $1$-cochain that cobounds $-\rho^* \chi$.  
From the cocycle condition (\ref{eq:homcond}) it follows that 
$f\colon G \to A$ is uniquely determined by its value on the generators. 
Now note that the map $\bar{f}=f \phi \colon F_n\to \G$ vanishes 
on the relators of $G$:
\begin{equation}
\label{eq:barf}
\bar{f}(r_k)=0,  \qquad 1\le k\le m.
\end{equation}
To finish the proof, we need to express this system of equations in terms 
of the the values $\bar{f}(x_i)=a_i$, and count the number of solutions. 

To that end, let $r=u_1 \cdots u_l $ be a word in $F_n$, with 
$u_j = x_{i_j}^{e_j}$. Then the following equality holds in the abelian group $A$:
\begin{equation}
\label{eq:foxlift}
\begin{split}
\bar{f}(r)=
\sum_{i=1}^{n} \s\bar\rho\big(\tfrac{\partial r}{\partial x_i}\big)(a_i) +
\sum_{j=1}^{l} & \tfrac{e_j-1}{2}\chi \big(
\bar\rho (u_{j} ), \bar\rho (u_{j}^{e_j} ) \big) +
\\
&\sum_{j=1}^{l-1}  \chi \big(
\bar\rho (u_1 \cdots  u_j ),
\bar\rho(u_{j+1} ) \big).
\end{split}
\end{equation}
This follows by induction on the length $l$ of the word $r$, 
using (\ref{eq:homcond}).  A detailed proof, in the case when 
$\Gamma=A\rtimes_{\s} B$, is given in \cite[Lemma~7.3]{MS}. 
The general case works similarly. 

From (\ref{eq:foxlift}), it is apparent that the system of equations 
(\ref{eq:barf}) coincides with (\ref{eq:foxnh}). By definition, the 
set of solutions to (\ref{eq:foxnh}) equals $Z^1_{\s\rho,\chi}(G,A)$.  
When this set is non-empty (i.e., $\epsilon_{\chi}(\rho)=1$), 
then $\abs{Z^1_{\s\rho,\chi}(G,A)}=\abs{Z^1_{\s\rho}(G,A)}$, 
cf.~\S\ref{subsec:counthom}.   So it is enough to count solutions of 
the system (\ref{eq:foxnh}) in the particular case when $[\chi]=0$.  
But this system is (\ref{eq:foxhom}), and we are done. 
\qed\end{pf}  

In particular, if $G=F_n$, then the system (\ref{eq:foxnh}) is 
empty, and so $\epsilon_{\chi}(\rho)=1$, for any homomorphism 
$\rho\colon F_n \to B$ and extension $\Gamma=A\times_{\s,\chi} B$. 

\subsection{Corank of twisted Jacobian}
\label{subsec:depth}

With notation as in Theorem~\ref{thm:fpcount}, suppose 
$A$ is an abelian $q$-group, where $q$ is a prime. 
Then the number of solutions of  system (\ref{eq:foxhom}) 
is of the form $q^d$, for some $d\ge 0$.  Denote this integer 
by $d(\s \rho)$.  Then:
\begin{equation}
\label{eq:z1rho}
\abs{Z^1_{\s\rho}(G, A)}= q^{d(\s\rho)} .
\end{equation}
For an arbitrary finite abelian group $A$, denote by $A_q$ the 
$q$-torsion subgroup.  
Then $A=\bigoplus_{q\mid\, \abs{A}} A_q$ as $G$-modules, 
and $Z^1_{\s\rho}(G, A)=\bigoplus_{q} Z^1_{\pi_q\s\rho}(G, A_q)$, where  
$\pi_q\colon \Aut(A)\to \Aut(A_q)$ is the canonical projection.  Hence: 
\begin{equation}
\label{eq:z1rhoq}
\abs{Z^1_{\s\rho}(G, A)}=\prod_{q\mid\, \abs{A}}  q^{d(\pi_q \s \rho)}.
\end{equation}

Now assume $A$ is homocyclic, say $A=\Z_{q^r}^{\oplus s}$. 
Then $\Aut(A)$ can be identified with $\GL(s, \Z_{q^r})$. 
Thus (\ref{eq:foxhom}) becomes a system of linear equations 
over the ring $\Z_{q^r}$, and 
\begin{equation}
\label{eq:drho}
d(\s\rho)=\corank \big( J_G^{\s\rho}\big) .
\end{equation}
Here recall $J_G=\phi \left(\partial r_i / \partial x_j\right)$ is 
the $m \times n$ matrix over $\Z G$ 
associated to the presentation $\mathcal{P}=\langle x_1,\dots, x_n \mid 
r_1,\dots, r_m \rangle$ for $G$, while $J_G^{\s\rho}$  is the $ms \times ns$ 
matrix over $\Z_{q^r}$ obtained by replacing 
each entry $e$ of $J_G$ by the matrix 
$\s\rho(e)\in \Aut(A)=\GL(s, \Z_{q^r})$. 

\section{Generalized Gasch\"{u}tz formula}
\label{sec:chief}

In this section, we give a formula counting the number of epimorphisms 
from a finitely presented group $G$ to a finite, 
solvable group $\G$.  We use \cite{Rob} as a general reference 
for group theory. 

\subsection{Complements in finite solvable groups}
\label{subsec:complements}

A group is said to be {\em solvable} if its derived series terminates. 
For a finite group $\G$, this is equivalent to $\G$ having 
an elementary abelian chief series, i.e., a non-refinable series 
of normal subgroups, such that  all the quotients are elementary 
abelian. By a classical result, any minimal normal subgroup of 
$\G$ must be an elementary abelian group $E$; moreover, the 
quotient $B=\G/E$ acts linearly on $E$.  

We will need the following result of Gasch\"{u}tz.  

\begin{theorem}[\cite{Ga}, Satz~3]
\label{thm:satz3}
Let $\G=E\times_{\sigma,\chi} B$ be an extension of a finite, 
solvable group $B$ by an  elementary abelian $q$-group  $E$ 
which is a minimal normal subgroup of $\G$.
Then the number of complements of $E$ in $\G$ is given by
\begin{equation}
\label{eq:cgamma}
c(\G)=c_{\chi}\cdot\abs{E}^{\zeta}\cdot q^{\kappa(\a-1)}, 
\end{equation}
where $c_{\chi}=1$ or $0$ according as $[\chi]=0$ or not,
$\zeta=0$ or $1$ according as $B$ acts
trivially on $E$ or not, and $\a$ is the number
of complemented chief factors of $\G$ isomorphic to $E$ as
$\G$-modules under the conjugation action.
\end{theorem}

If $c_{\chi}=0$, the extension is non-split, and so $c(\G)=0$; 
the case $c_{\chi}=1$ is the one requiring an argument.  
The proof in \cite{Ga} breaks into several steps.
First, it is shown that $c(\G)=c(\G/Z)\cdot c(Z)$,
where $Z=C_{\G}(E)$ is the centralizer of $E$ in $\G$,
$c(\G/Z)$ is the number of complements of $E$ in $\G/Z$,
and $c(Z)$ is the number of complements of $E$ in $Z$.
Next, it is shown that $c(\G/Z)=\abs{E}^{\zeta}$, and
$c(Z)=\abs{\End_{\G}(E)}^{\a-1}$.  Finally, it is noted that 
$\abs{\End_{\G}(E)}=q^{\kappa}$, for some $\kappa\ge 0$.

For related results, see \cite[(2.10)]{AG} and 
\cite[Theorem~2]{DL}.

\subsection{A recursion formula}
Now let $G$ be a finitely presented group. 
\begin{lemma}
\label{lem:key}
Suppose $\G=E\times_{\sigma,\chi} B$ is an extension of a finite group
$B$ by an elementary abelian $q$-group $E$ which is a minimal 
normal subgroup of $\G$. Then: 
\begin{equation}
\label{eq:indstep}
\abs{\Epi(G,\G)}=\sum_{\rho\in\Epi(G,B)}
 \left(  \epsilon_{\chi}(\rho) q^{d(\s \rho)} - c\right), 
\end{equation}
where $c=c(\G)=c_{\chi} \abs{Z^1_{\s}(B,E)}$ is the 
number of complements of $E$ in $\G$ ($c_{\chi}=1$ 
if the extension splits, in which case $\epsilon_{\chi}(\rho)=1$, and 
$c_{\chi}=0$ otherwise).
\end{lemma}

\begin{pf}
Fix an epimorphism $\rho\colon G\to B$.  
Then $\rho$ has  $\abs{Z_{\s\rho,\chi}^{1}(G, E)} =
\epsilon_{\chi}(\rho) \cdot $ $ q^{d(\s\rho)}$ lifts to $\G$.  
Let  $\l\colon G\to \G$ be such a lift, and let $U=\Im\l$. Then 
$U$ is an extension of $B$ by $K=U\cap E$ (a subgroup of $E$). 
By minimality, $B$ acts irreducibly on $E$, and so either $K=E$,  
in which case  $U=\G$ (and so $\l$ is surjective), or $K=1$, 
in which case $U$ is a complement of $E$.  
Therefore $\rho$ contributes $\epsilon_{\chi}(\rho) q^{d(\s\rho)}-c$  
to $\abs{\Epi(G,\G)}$. 
\qed\end{pf}

We now have:
\begin{equation}
\label{eq:qdrho}
\begin{split}
q^{d(\s\rho)} &= \abs{Z^1_{\s\rho} (G,E)} \\
 &=  \abs{B^1_{\s\rho} (G,E)} \cdot \abs{H^1_{\s\rho} (G,E)} \\
& =  \abs{E}^{\zeta} \cdot q^{\dim_{\Z_q} H^1_{\s\rho} (G,E)}. 
 \end{split}
\end{equation}

Combining Lemma~\ref{lem:key} with Theorem~\ref{thm:satz3} 
and formula (\ref{eq:qdrho}), we obtain the following. 

\begin{theorem}
\label{thm:gg}
Let $G$ be a finitely presented group, and
let $\G=E\times_{\sigma,\chi} B$ be an extension of a finite, solvable 
group $B$ by an  elementary abelian $q$-group  $E$ which is  a 
minimal normal subgroup of $\G$. Then:
\[
\abs{\Epi(G,\G)}=\abs{E}^{\zeta}  \sum_{\rho\in\Epi(G,B)}
\left( \epsilon_{\chi} (\rho) q^{\dim_{\Z_q}H^{1}_{\s\rho}(G,E)} - 
c_{\chi} q^{\kappa(\a-1)} \right).
\]
\end{theorem}

This Theorem allows us to compute recursively 
$\abs{\Aut(\G)}=\abs{\Epi(\G,\G)}$, 
for any finite solvable group $\G$. 

\begin{corollary}
\label{cor:aut}
Let $\G=E\times_{\sigma,\chi} B$ be an extension of a finite, solvable 
group $B$ by an  elementary abelian $q$-group  $E$ which is  a minimal 
normal subgroup of $\G$. Then:
\[
\abs{\Aut(\G)}=\abs{E}^{\zeta}  \sum_{\rho\in\Epi(\G,B)}
\left( \epsilon_{\chi} (\rho) q^{\dim_{\Z_q}H^{1}_{\s\rho}(\G,E)} - 
c_{\chi} q^{\kappa(\a-1)} \right).
\]
\end{corollary}

Combining the two results above gives a recursion formula 
for the Hall invariant $\delta_{\G}(G)$, for any finite solvable group $\G$. 

\subsection{Lifting through the chief series} 
\label{subsec: lift chief}
We now describe an explicit procedure for constructing 
the set $\Epi(G,\G)$, and counting its elements. 
Start with a chief series 
\begin{equation}
\label{eq:chiefseries}
\G=\G_0>\G_{1}> \cdots> \G_{\nu}> \G_{\nu+1}=1. 
\end{equation}
Write $E_i=\G_i/\G_{i+1}=\Z_{q_i}^{\oplus s_i}$ and $B_i=\G/\G_{i}$, 
for $0\le i\le \nu$. 
Let $\chi_i\colon B_i\times B_i \to E_i$ be a 
classifying $2$-cocycle for the extension 
\begin{equation}
\label{eq:bext}
\xymatrix{
1\ar[r] & E_i\ar[r] & B_{i+1} \ar[r] &  B_i \ar[r] &  1,}
\end{equation}
with monodromy $\s_i\colon B_i \to \Aut(E_i)=\GL(s_i,q_i)$.  
Finally, let  $c_i=c(B_{i+1})$ be the number of complements 
of $E_i$ in $B_{i+1}$, let $\a_i$ be the number of chief factors 
of $B_{i+1}$ isomorphic to $E_i$ as $B_{i+1}$-modules, and 
set $q_i^{\kappa_i}=\abs{\End_{B_{i+1}} (E_i)}$. 

Now let $G$ be a finitely presented group. Start with 
the trivial epimorphism $\rho_0\colon G \to B_0=1$. 
Then $\Epi(G,B_1)$ consists of $q_0^{\beta_0}-1$ elements, 
where $\beta_0=\dim_{\Z_{q_0}} H^1(G; E_0)$.   

For each such element $\rho_1\colon G\surj B_1$, we must 
decide whether there is a lift $\rho_2\colon G\surj B_2$. 
Applying Theorem~\ref{thm:gg} to the extension 
$B_2=E_1\times_{\s_1,\chi_1} B_1$, we find there 
are $\abs{E_1}^{\zeta_1}  \sum_{\rho_1\in\Epi(G,B_1)}
\big( \epsilon_{\chi_1} (\rho_1) q_1^{\beta_1} - 
c_{\chi_1} q_1^{\kappa_1(\a_1-1)} \big)$ such lifts, where 
$\beta_1=\dim_{\Z_{q_1}}H^{1}_{\s_1\rho_1}(G,E_1)$. 
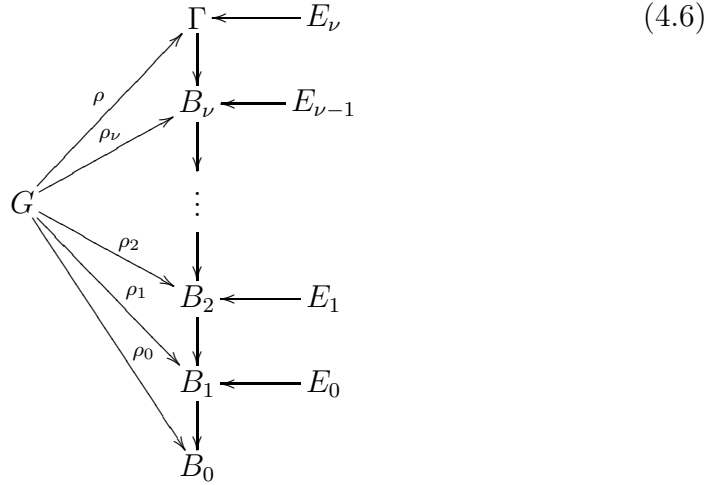
\begin{equation}
\label{eq:lift diagram}
\xymatrixrowsep{18pt}
\xymatrixcolsep{22pt}
\xymatrix{
&&  \G  \ar[d]  & E_{\nu} \ar[l]\\
 &&  B_{\nu} \ar[d]  & E_{\nu-1} \ar[l]\\
G   \ar[rruu]^(.5){\rho} \ar[rru]^(.56){\rho_{\nu}} 
\ar[rrd]^(.56){\rho_2} \ar[rrdd]^(.56){\rho_1} \ar[rrddd]^(.6){\rho_0}
 &&\! \:^{\vdots} \ar[d]  \\
&& B_2 \ar[d] & E_{1} \ar[l] \\
& & B_1\ar[d] & E_{0} \ar[l] \\
&& \: B_0 \:
}
\end{equation}

Continuing in the manner illustrated in diagram (\ref{eq:lift diagram}), 
we obtain the following formula.

\begin{corollary}
\label{cor:gengas}
With notation as above, the number of epimorphisms 
$\rho\colon G \surj \G$ is given by
\[
\abs{\Epi(G,\G)}=
\sum_{\rho_{1}\in\Epi_{\rho_{0}}(G,B_{1})} \!\cdots\!
\sum_{\rho_{\nu}\in\Epi_{\rho_{\nu-1}}(G,B_{\nu})}
\abs{E_{\nu}}^{\zeta_{\nu}} \! \left( \epsilon_{\chi_{\nu}}(\rho_{\nu})  
q_{\nu}^{\beta_{\nu}} - 
c_{\chi_{\nu}} q_{\nu}^{\kappa_{\nu}(\a_{\nu}-1)} \right) \! ,
\]
where  $B_{\nu+1}=E_{\nu}\times _{\s_\nu, \chi_\nu} B_{\nu}$, 
$\beta_{\nu}=\dim_{\Z_{q_{\nu}}} 
H^{1}_{\s_{\nu}\rho_{\nu}} (G, E_{\nu})$, 
$\alpha_{\nu}$ is the number of complemented chief 
factors in $B_{\nu+1}$, isomorphic to $E_{\nu}$ as 
$B_{\nu+1}$-modules, and $\Epi_{\rho_{i}}(G,B_{i+1})$ is 
the set of epimorphisms lifting $\rho_{i}\colon G\surj B_{i}$.
\end{corollary}

When $G$ is the free group of rank $n$, Theorem~\ref{thm:gg} 
(or Corollary \ref{cor:gengas}) recovers the classical Gasch\"{u}tz formula.   

\begin{corollary}[\cite{Ga}, Satz 4]
The Eulerian function of a finite solvable \linebreak group, 
$\phi(\Gamma,n)=\abs{\Epi(F_n,\Gamma)}$, is given by 
\[
\phi(\Gamma,n)=\prod_{i=1}^h \left[ q_i^{s_i v_i n} 
\left(q_i^{s_i n}-q_i^{s_i \zeta_i}\right)\left(q_i^{s_i
n}-q_i^{s_i\zeta_i+\kappa_i}\right)\cdots \left(q_i^{s_i
n}-q_i^{s_i\zeta_i+(u_i-1)\kappa_i}\right)\right], 
\]
where $V_1, V_2,\dots, V_h$ are the distinct $\Gamma$-module 
isomorphism types of the chief factors, $q_i^{s_i}=\abs{V_i}$,
$\zeta_i=0$ or $1$ according as $\G$ acts trivially on $V_i$ or not,
$q_i^{\kappa_i}$ is the number of  $\G$-endomorphisms of $V_i$, 
and finally, $u_i$ is the numbers of factors of type $V_i$ 
which are complemented in $\G$, and $v_i$ is the number 
of  other factors of type $V_i$.
\end{corollary}

\begin{pf}
Let $E$ be a minimal normal subgroup of $\G$, and set $B=\G/E$. 
Let $V_i$ be the $\G$-isomorphism type of $E$.  
Applying Theorem~\ref{thm:gg} we find an expression of the form:
\[
\abs{\Epi(F_n,\G)}=\abs{\Epi(F_n,B)} \cdot \abs{E}^{\zeta_i} 
\big(  q_i^{\beta_i}-c_{\chi_i} q_i^{\kappa_i (\alpha_i-1)} \big) , 
\]
where $\beta_i=\dim_{\Z_{q_i}} H^1(F_n, E)=s_i(n-\zeta_i)$.  
If $E$ is not complemented, then $c_{\chi_i}=0$, and 
the second factor reduces to 
$q_i^{s_i n}$; otherwise, the second factor reduces to 
$q_i^{s_i n} - q_i^{s_i \zeta_i+(u_i-1)\kappa_i}$. 

Now repeat the procedure, with $\G$ replaced by $B$, and 
continue in this fashion till the trivial group is reached.  
\qed\end{pf}

\section{Dihedral groups}
\label{sec:dihedral}

We now study in more detail epimorphisms to  
the dihedral group of order $2m$.  This group is a split 
extension $D_{2m}=\Z_m\rtimes_{\s} \Z_2$, with monodromy 
$\sigma(b)=b^{-1}$:
\[
D_{2m}=\langle a,b \mid a^m=b^2=1, \, bab=a^{-1}\rangle.
\]
Let $m=q_1^{\a_1}\cdots q_r^{\a_r}$ be the prime decomposition of $m$.  
A chief series for $D_{2m}$ is
\[
D_{2m} =\G_0 > \G_{1,1} > \cdots > \G_{1,\a_1} > \cdots 
> \G_{r,1} > \cdots > \G_{r,\a_r} >  1,
\]
with terms $\G_{i,j} = \Z_{m/( q_1^{\a_1} \cdots 
q_{i-1}^{\a_{i-1}}q_i^{j-1})}$.  The chief factors are 
$E_0=\Z_2$ and $E_{i,j}=\Z_{q_{i}}$.  
The lifting process goes through the extensions $B_1=\Z_2$ and 
$B_{i,j} = D_{2 q_1^{\a_1} \cdots q_{i-1}^{\a_{i-1}}q_i^{j-1}}$, 
for $1\le i \le r$, $1\le j\le \a_i$.  Of these extensions, only the 
ones where a prime $q_i$ appears for the first time are split, 
while the others are non-split.  Indeed, all the extensions are 
of the form $D_{2ql}=\Z_q\rtimes_{\s,\chi} D_{2l}$, 
with $\chi (a^u b^v, a^s b^t)=k$, 
where $u+s\cdot (-1)^v= l\cdot k + r \md{ql}$,
and $0\le r < l$.  

As before, let $G$ be a finitely presented group. 
Applying Lemma \ref{lem:key}, we obtain the following recursion formula 
for the number of epimorphisms from $G$ to a dihedral group:
\begin{equation}
\label{eq:epi dihedral}
\abs{\Epi(G,D_{2ql})}=
\left\{
\begin{array}{ll}
\sum_{\rho\in\Epi(G,D_{2l})} \big( q^{d(\s\rho)}-q\big) &
\hbox{if $q\nmid l$},
\\ 
\sum_{\rho\in\Epi(G,D_{2l})} \epsilon_{\chi}(\rho) q^{d(\s\rho)} &
\hbox{if $q \mid l$} .
\end{array}
\right.
\end{equation}

\begin{example}
\label{ex:free dihedral}
{\rm
For the free group $G=F_n$, we find:  
\begin{equation}
\label{eq:phi dihedral}
\abs{\Epi(F_n,D_{2m})} = (2^n-1) m^n \cdot 
\prod_{i=1}^{r}  (1- q_i^{1-n}),
\end{equation}
where $q_1,\dots , q_r$ are the prime factors of $m$ 
(and the empty product is  $1$). 
This recovers a computation of Kwak, Chun, 
and Lee, see \cite[Lemma~4.1]{KCL}.
}
\end{example}

\begin{example}
\label{ex:surface}
{\rm
Let $\Pi_g=\langle x_1,\dots ,x_g,y_1,\dots, y_g \mid 
[x_1,y_1]\cdots [x_g,y_g]=1\rangle$ and 
$\Pi_g^*=\langle x_1,\dots ,x_g \mid 
x_1^2\cdots  x_g^2=1\rangle$ be the fundamental groups 
of orientable (respectively, non-orientable) surfaces of genus $g$. 
Let $q_1,\dots, q_r$ be the odd prime factors of $m$, and put 
$e=m/2 \md{2}$ if $m$ is even.  
Then, according to whether $m$ is odd or even, 
\begin{gather}
\label{eq:mg dihedral}
\abs{\Epi(\Pi_g,D_{2m})} =
\begin{cases}
\DS{ m^{2g-1}  (2^{2g}-1) \prod_{i=1}^{r}   \big(1- q_i^{2-2g}\big),}
 \\[3pt]
\DS{m^{2g-1}   (2^{2g}-1)  (2^{e}- 2^{2-2g}) \prod_{i=1}^{r}  
 \big(1- q_i^{2-2g}\big),}
\end{cases}
\\[6pt]
\nonumber
\abs{\Epi(\Pi_g^*,D_{2m})} =
\begin{cases}
\DS{ m^{g-1}  \Big[(2^{g}-2)   \prod_{i=1}^{r}  
\big(1- q_i^{2-g}\big)+ \prod_{i=1}^{r}   \big(q_i- q_i^{2-g})\Big],}
  \\[3pt]
 \DS{  m^{g-1} (2^{e}- 2^{2-g}) 
 \Big[(2^{g}-2)   \prod_{i=1}^{r} \big(1- q_i^{2-g}\big)+ \prod_{i=1}^{r}
 \big(q_i- q_i^{2-g}\big)\Big].}
\end{cases}
\end{gather}

In \cite{KL}, Kwak and Lee obtained 
related formulas, counting the number of regular, $D_{2p}$ 
branched covers ($p$ prime)  of a closed surface. 
}
\end{example}

\begin{example} %Baumslag-Solitar groups}
\label{sec:baumslag solitar d8}
{\rm
An interesting family of examples is provided by 
the Baumslag-Solitar one-relator groups 
$\BS(m,n)=\langle x,y \mid x y^{m} x^{-1} y^{-n} \rangle$. 
As an illustration of our techniques, 
we now compute the number of epimorphisms from 
$G=\BS(m,n)$ to $D_8$.   

The dihedral group $D_8$ is a central extension of 
$\Z_2^{\oplus 2}$ by $\Z_2$, with $2$-cocycle 
$\chi\colon \Z_2^{\oplus 2} \times \Z_2^{\oplus 2} \to \Z_2$ 
assuming non-zero values only on the pairs 
$(a,a)$, $(b,a)$, $(a,ab)$, and $(b,ab)$.  An epimorphism 
$G\surj D_8$ induces by abelianization an epimorphism
$\Z\oplus\Z_{\abs{n-m}}\surj \Z_2\oplus\Z_2$; this can happen 
only if $m$ and $n$ have the same parity.

So assume $m  \equiv n\md{2}$.  Then there are 
precisely $6$ epimorphisms from $G$ to  
$\Z_2^{\oplus 2}$; let $\rho\colon G\surj \Z_2^{\oplus 2}$ 
be one of those.  
A computation shows that $J_G^{\rho}=0$.  
By  Theorem \ref{thm:fpcount},   
$\rho$ lifts to $D_8$, i.e., $\epsilon_{\chi}(\rho)=1$, 
if and only if 
\[
\sum_{k=1}^{m} \chi(u v^{k-1}, v)-\chi(u, u)+\chi(u v^{m},u)+
\sum_{l=1}^{n}\big( -\chi(v, v)+\chi(u v^m u v^{l-1} ,v)\big) = 0.
\]
where $u=\rho(x)$ and $v=\rho(y)$, in which case there are 
precisely $4$ lifts to $D_8$.   
This equation simplifies to 
\begin{align*}
&\tfrac{m}{2}\chi(u,v)+\tfrac{m}{2}\chi(uv,v)+\tfrac{n}{2}\chi(v,v)= 0, 
\quad \hbox{or} 
\\[3pt]
&\tfrac{m+1}{2}\chi(u,v)+\tfrac{m-1}{2}\chi(uv,v)+\tfrac{n-1}{2}\chi(v,v)
+\chi(uv,u)-\chi(u, u)= 0,
\end{align*}
according to whether $m$ is even or odd. 
We find that $\abs{\{ \rho \mid \epsilon_{\chi} (\rho)=1\}} = 
6, 4, 2,\:  \hbox{or } 0$, according to whether $(m, \tfrac{n-m}{2})
\equiv (0,0), (0,1), (1,1),\: \hbox{or } (1,0)$ modulo $2$, respectively. 
Using the fact that $\Aut (D_8) =D_8$, we obtain:
\begin{equation}
\label{eq:baumsol d8}
\delta_{D_8}( \BS(m,n) ) = 
\begin{cases}
3  &\hbox{if $m$ even and $n-m\equiv 0 \md{4} $,}\\
2  &\hbox{if $m$ even and $n-m\equiv 2 \md{4} $,}\\
1  &\hbox{if $m$ odd\ \ and $n-m\equiv 2 \md{4} $,}\\
0  &\hbox{otherwise}. 
\end{cases}
\end{equation}
}
\end{example}

\section{Binary dihedral groups}
\label{sec:bin dihedral}

The binary dihedral group of order $4m$ is a central extension 
$D^*_{4m}=\Z_2\times _{\chi} D_{2m}$, with presentation
\[
D^*_{4m}=\langle a, b \mid a^{2m}=1, \, a^m=b^2,\,  
bab^{-1}=a^{-1}\rangle. 
\]
In particular, $D^*_4=\Z_4$, $D^*_8=Q_8$, the quaternion group, and 
$D^*_{4\cdot 2^{m-2}}=Q^{}_{2^m}$, the generalized quaternion group. 
Write $m=q_0^{\a_0}q_1^{\a_1}\cdots q_r^{\a_r}$, with $q_0=2$.  
A chief series is then:
\[
D^*_{4m}=\G_0 > \G_{0,0}  > \G_{0,1} > \cdots > \G_{0,\a_0}  > \cdots 
> \G_{r,1} > \cdots > \G_{r,\a_r} >  1,
\]
with terms $\G_{0,0} = \Z_{2m}$, 
$\G_{i,j} = \Z_{m/( q_0^{\a_0} \cdots q_{i-1}^{\a_{i-1}}q_i^{j-1})}$ 
and factors $E_{0,0}=\Z_2$, $E_{i,j}=\Z_{q_{i}}$, 
where $0\le i \le r$, $1\le j\le \a_i$. 
The lifting process goes through the extensions
$B_0=\Z_2$, $B_{0,j}=D_{2^{j+2}}$, and 
$B_{i,j} = D^*_{4 q_0^{\a_0} \cdots q_{i-1}^{\a_{i-1}}q_i^{j-1}}$. 
Here only the extensions where a prime $q_i$ appears for 
the first time are split, the rest are non-split, except 
when $m$ is even, in which case the prime $q_0=2$ produces 
a split extension the first two times it appears.  

Indeed, there are three types of extensions that occur:
\begin{itemize}
\item $D^*_{2ql}  = \Z_q \times _{\s, \chi} D^*_{2l}$, with 
$q$ an odd prime, in which case the computation 
of $\chi$ goes essentially as in the dihedral case.   

\item $D_{2^{r+1}} = \Z_2 \times_{\chi}  D_{2^{r}}$, for which 
$\chi$ was computed before.

\item  $Q_{2^{r+1}} = \Z_2 \times _{\chi} D_{2^{r}}$, 
in which case $\chi (a^u b^v, a^s b^t)=k+l$, where
$u+s\cdot (-1)^v \equiv k\cdot 2^{r-1} + n \md{2^r}$ 
with $0\le n<  2^{r-1}$, and $l=1$ if $v=t=1$ and $l=0$ otherwise.  

\end{itemize}

Applying Lemma \ref{lem:key}, we obtain the following recursion formulas:
\begin{align}
\abs{\Epi(G,D^*_{2ql})} &=
\begin{cases}
\DS{\sum_{\rho\in\Epi(G,D^*_{2l})} \big( q^{d(\s\rho)}-q\big)} &
\hbox{if $q\nmid l$},
\\[6pt]
\DS{\sum_{\rho\in\Epi(G,D^*_{2l})} \epsilon_{\chi}(\rho) q^{d(\s\rho)} }&
\hbox{if $q \mid l$} ,
\end{cases}
\\[6pt]
\abs{\Epi(G,Q_{2^{r+1}})} &=
\sum_{\rho\in\Epi(G,D_{2^r})} \epsilon_{\chi}(\rho) q^{d(\rho)}.
\end{align}

\begin{example}
\label{ex:free binary dihedral}
{\rm
For the free group $G=F_n$,  we find: 
\begin{equation}
\label{eq:phi binary dihedral}
\abs{\Epi(F_n,D^*_{4m})} = (4^n-2^n) m^n \cdot 
\prod_{i=0}^{r}  (1- q_i^{1-n}) . 
\end{equation}
}
\end{example}

\begin{example} %Baumslag-Solitar groups}
\label{sec:baumslag solitar}
{\rm 
Let us compute the number of epimorphisms from the 
Baum\-slag-Solitar groups to the quaternion group. 
Notice that $Q_8=\Z_2\times_{\chi} \Z_2^{\oplus 2}$, with 
$2$-cocycle $\chi$ vanishing only on the pairs 
$(a,b)$, $(b,ab)$, and $(ab,a)$.
Let $\rho\colon \BS(m,n)\surj \Z_2\oplus \Z_2$ be 
an epimorphism.  Then necessarily $m$ and $n$ have the 
same parity. Moreover, 
$\rho$ lifts to $Q_8$ if and only if
\[
\begin{split}
\sum_{k=1}^{m} \chi(u v^{k-1}, v)&-\chi(u, u)+\chi(u v^{m},u)+
\\ 
&\sum_{l=1}^{n}\big( -\chi(v, v)+
\chi(u v^m u v^{l-1} ,v) \big) =0, 
\end{split}
\]
where $u=\rho(x)$ and $v=\rho(y)$, in which case  
$\rho$ has $4$ lifts.  
The above condition is equivalent to $m+n \equiv 0\md{4}$. 
Using the fact that $\Aut(Q_8)=S_4$, we conclude:
\begin{equation}
\label{eq:qbaumsol}
\delta_{Q_8}( \BS(m,n) ) = 
\begin{cases}
1  &\text{if $n-m$ is even and $m+n\equiv 0\md{4}$},\\
0   &\text{otherwise}. 
\end{cases}
\end{equation}
}
\end{example}

\section{Finite quotients of derived length $3$}
\label{sec:sym4}

We now consider epimorphisms of a finitely presented group $G$ onto
finite groups which are not metabelian. A nice class of groups of 
derived length $3$ are the split extensions $\G=\Z_q^2\rtimes_{\s} D_{2p}$, 
where $p$ and $q$ are distinct primes such that $q$ has order $2 \md{p}$.
On generators $b$ and $c$ for $D_{2p}$, 
the monodromy $\s\colon D_{2p}\to \GL(2,q)$ is given by 
$\s(b)=\big( \begin{smallmatrix} r & 1 \\ -1 & 0
 \end{smallmatrix}\big)$
and $ \s(c)=\big( \begin{smallmatrix}
 0 & 1 \\ 1 & 0 \end{smallmatrix}\big) $, for some $r$.
  
According to Theorem \ref{thm:gg}, we have:
\begin{equation}
\label{eq:epipq}
\abs{\Epi(G,\G)}=  q^2\sum_{\rho\in\Epi(G,D_{2p})}
\big(q^{\beta(\rho)} - 1 \big), 
\end{equation}
where $\beta(\rho)=\dim_{\Z_q}H^{1}_{\sigma\rho}(G;\Z_q^{\oplus 2})$.  

In particular, the symmetric group on four letters is a split extension 
$S_4=\Z_2^{\oplus 2} \rtimes_{\s} S_3$, 
with monodromy $\sigma\colon S_3=\SL(2,2)\to 
\Aut(\Z_2^{\oplus 2})=\GL(2,2)$ the natural inclusion, given by 
$\s(b)=\big( \begin{smallmatrix} 1 & 1 \\ 1 & 0
 \end{smallmatrix}\big)$
and $ \s(c)=\big( \begin{smallmatrix}
 0 & 1 \\ 1 & 0 \end{smallmatrix}\big) $.  We then have:
$\abs{\Epi(G,S_4)}=  4\sum_{\rho\in\Epi(G,S_3)}
\big(2^{\beta(\rho)} - 1 \big)$,   
where $\beta(\rho)=\dim_{\Z_2}H^{1}_{\sigma\rho}(G;\Z_2^{\oplus 2})$.  
Since $\Aut(S_4)=S_4$, we find:
\begin{equation}
\label{eq:del s4}
\delta_{S_4}(G)=\tfrac{1}{6} 
\sum_{\rho\in\Epi(G,S_3)}\big(2^{\beta(\rho)} - 1 \big).
\end{equation}

\begin{example}
\label{ex:epi free s4}
{\rm
Consider $G=F_n$.  Recall that
$\abs{\Epi(F_n,S_3)}=(2^n-1)(3^{n}-3)$. If $\rho\colon F_n \surj S_3$, 
then $H^{1}_{\sigma\rho}(F_n;  \Z_2^{\oplus 2})=\Z_2^{\oplus 2n-2}$. 
Hence, 
\begin{equation}
\label{eq:eulers4}
\abs{\Epi(F_n,S_4)}=(2^n-1)(3^{n}-3) (4^{n} -4).
\end{equation}
}
\end{example}

\begin{example} %[Surface groups]
\label{ex:surfaces}
{\rm
Let $G=\langle x, y\mid yxy^{-1} =x^{-1}\rangle$ 
be the Klein bottle group.  There is then an obvious epimorphism 
$\rho\colon G \surj S_3$; in fact, $\delta_{S_3}(G)=1$.  But this 
epimorphism does not lift to $S_4$.  Indeed, 
$J_G=( \begin{array}{cc}  y+x^{-1} &\: 1-x^{-1}\end{array} )$, 
and thus $J_{G}^{\s\rho}=\big( \begin{smallmatrix}  1&0&0&1 \\ 
0&0&1&1 \end{smallmatrix} \big)$. 
Hence, $H^1_{\s\rho}(G,\Z_2^{\oplus 2})=0$, and 
so $\delta_{S_4}(G)=0$. 
}
\end{example}

\begin{example} %{Baumslag groups}
\label{ex:baumslag}
{\rm
We now show that the Baumslag parafree groups, 
$\P(m,n)=\langle x,y,z \mid x z^{m} x z^{-m} x^{-1} z^{n} y z^{-n} y^{-1} \rangle$, 
fall into at least two distinct isomorphism classes, 
each containing infinitely many members. We do this 
by counting epimorphisms to $S_4$. 

Let us start by computing the Fox Jacobian of $G=\P(m,n)$:
\[
J_G=
\begin{pmatrix}
1+xz^{m}-[y, z^{n}]  \
&  yz^{n}y^{-1}-1  \
& x(1-z^{m} x z^{-m})s_m+([y, z^{n}]-y)s_n
\end{pmatrix},
\]
where $s_k=1+z+\cdots+z^{k-1}$.
The abelianization $\ab\colon G\to G/G'=\Z^2=\langle t_1, t_2\rangle$ 
sends $x\mapsto 1$, $y\mapsto t_1$, $z\mapsto t_2$. 
Thus,  
\[
J_G^{\ab}=\begin{pmatrix}
t_2^{m}  \quad
&  t_2^{n}-1 \quad
& (1-t_1)(1+t_2+\dots+t_2^{n-1})
\end{pmatrix}.
\]

Clearly, $\abs{\Epi(G,\Z_2)}=3$. Since the 
first entry in $J_G^{\ab}$ is never zero, each epimorphism
$ G \surj \Z_2$ lifts to $6$ different epimorphisms to $S_3$. 
Writing the typical element of $S_3=\Z_3\rtimes \Z_2$ as 
$b^{\b}c^{\g}$, we see that the $18$ epimorphisms 
$\rho\colon G\surj S_3$ divide into $3$ families:
\begin{itemize}
\item $x\to b^\b, y\to b^{\b_1} c, z\to b^{\b_2}$, 
where $\b=n\b_2\ne 0$ and $\b_1$ arbitrary, 

\item $x\to b^\b, y\to b^{\b_1}, z\to b^{\b_2} c$, 
where $\b=((-1)^m-(-1)^{m+n})\b_1\ne 0$ and $\b_2$ arbitrary, 

\item $x\to b^\b, y\to b^{\b_1} c, z\to b^{\b_2}$, 
where $\b=((-1)^m-(-1)^{m+n})(\b_1-\b_2)\ne 0$. 
\end{itemize}
Each $\rho$ in the first family contributes $4(2^2-1)$ to 
$\abs{\Epi(G,S_4)}$, as $J^{\s\rho}$ has corank $4$, while 
the other $\rho$'s contribute either $4(2^3-1)$ or $4(2^2-1)$, 
according as $J^{\s\rho}$ has corank $5$ or $4$ 
(depending  on whether $m-n \equiv 2\md{4}$ or not).  
Therefore
\begin{equation}
\label{eq:delbaum}
\delta_{S_4}( \P(m,n) ) = 
\begin{cases}
\frac{6\cdot 4(2^2-1)+12 \cdot 4(2^3-1)}{24}=17 
  &\text{if $m$ odd, $m-n \equiv 2\md{4}$},\\[3pt]
 \frac{18\cdot 4(2^2-1)}{24}=9 &\text{otherwise}. 
\end{cases}
\end{equation}

It would be interesting to see whether solvable 
Hall invariants completely classify the Baumslag parafree 
groups, and, more generally, the parafree groups considered 
by Strebel in \cite{Str}. 
}
\end{example}

\begin{example} %[Link groups]
\label{ex:links}
{\rm 
Let $L$ be the $2$-component link from 
\cite[p.~72]{Hil}, and $G$ the fundamental group
of its complement, with presentation
\[
G=\left\langle x_1,x_2,x_3,x_4 \left|
\begin{array}{l}
 x_1^{x_4^{-1}} (x_4 x_1)^{x_2^{-1}x_1 x_2},\quad 
 (x_3^{-1} x_1 x_4)^{x_2^{-1} x_1^{-1}x_2} (x_1 x_4)^{x_3} \\
\left[x_1^{-1}x_4^{-1}x_3x_1x_4x_3^{-2}x_4, x_2\right] 
\end{array}\right.\right\rangle,
\]
where $x^y =y^{-1}xy$.  
Let $G_{\omega}$ be the intersection of the lower central series of $G$.  
As noted by Hillman, $G=G_{\omega} \rtimes \P(-1,1)$.  In particular, 
$L$ is not a homology boundary link.  Moreover, 
$G/G''\cong F_2/F_2''$, yet $G\not\cong F_2$.  

Notice that $\abs{\Epi(G,\G)}=\abs{\Epi(F_2,\G)}$, 
for any finite metabelian group $\Gamma$. In particular, 
$\delta_{S_3}(G)=3$. On the other hand, we can distinguish 
$G$ from both $F_2$ and $\P(-1,1)$ by counting 
representations onto $S_4$:
\[
\delta_{S_4}(G)= 2\cdot (2^4 -1)+ (2^2-1)=33. 
\]
}
\end{example}

\section{The lattice of subgroups}
\label{sec:euler}

Let $\G$ be a finite group.  Let $L(\G)$ be the lattice of 
subgroups of $\G$, ordered by inclusion.  The 
{\em M\"{o}bius function}, $\mu\colon L(\G) \times 
L(\G)\to \Z$, is defined inductively by $\mu(H,H)=1$, 
and $\sum_{H\le S \le K} \mu(H,S)=0$, for any subgroup 
$K\le \Gamma$.  For simplicity, write $\mu(H):=\mu(H,\G)$.  

In \cite{KT}, Kratzer and Th\'evenaz give a formula for the 
M\"obius function of a solvable group, in terms of a chief series. 

\begin{theorem}[\cite{KT}]
\label{thm:KT}
Suppose $\G$ is solvable, and 
let $\G=\G_0 > \cdots > \G_{\nu} > 1$ be a chief series.   
If $H\le \G$, let $H_i=\G_i H$, 
and consider the sequence $H=H_r<\cdots<H_0=\G$, 
where one keeps only distinct terms $H_i$. 
Let $h_i$ be the number of complements of $H_i$ in $L(\G)$ which 
contain $H_{i+1}$. Then 
\[
\mu(H)=(-1)^r  h_1\cdots h_{r-1}.
\]
\end{theorem}

In the particular case when $\G$ is nilpotent, this recovers 
a classical result of Weisner:  $\mu(H)=0$, unless 
$H\triangleleft \G$ and $\G/ H \cong 
\bigoplus_{i=1}^{r} \Z_{q_i}^{\oplus s_i}$, in which case 
$\mu(H)=\prod _{i=1}^{r} (-1)^{s_i} q_i^{s_i(s_i-1)/2}$. 

Now let $G$ be a finitely generated group.  
Then, as noted by P.~Hall \cite{HaP}, 
\begin{equation}
\label{eq:hallenum1}
\abs{\Hom(G,\G)}=\sum_{H\le\G} \abs{\Epi(G,H)},
\end{equation}
or, by M\"obius inversion:
\begin{equation}
\label{eq:hallenum2}
\abs{\Epi(G,\G)}=\sum_{H\le\G}\mu(H) \abs{\Hom(G,H)}.
\end{equation}

The {\em Eulerian function} of $\G$ is the sequence 
$\phi(\G,n)=\abs{\Epi(F_n,\G)}$, counting ordered $n$-tuples 
generating $\G$.  By the Hall enumeration principle 
(\ref{eq:hallenum2}), the Eulerian function is determined 
by the M\"obius function, as follows: 
$\phi(\G,n)=\sum_{H\le \G} \mu(H) \abs{H}^n$.  

In conjunction with the results from Section \ref{sec:homology}, 
the Hall enumeration principle provides 
an alternate way to compute the number of epimorphisms 
from an arbitrary finitely presented group $G$ to a finite 
solvable group $\G$.  

\begin{example}
\label{ex:lattice dihedral}
{\rm 
The lattice of subgroups of the dihedral group $D_{2m}$ 
consists of one subgroup of type $\Z_l$ and $m/l$ subgroups 
of type $D_{2l}$, for each divisor $l$ of $m$.  The M\"obius 
function is given by
\[
\mu(\Z_l)=-\tfrac{m}{l}\mu(\tfrac{m}{l}),
\qquad
\mu(D_{2l})=\mu(\tfrac{m}{l}).
\]
Let $q_1,\dots ,q_r$ be the prime divisors of $m$. 
By Proposition \ref{prop:homcount} and  formula 
(\ref{eq:z1rho}), we have:
\begin{eqnarray}
\label{eq:epid2m}
\nonumber
\abs{\Epi(G,D_{2m})} & = \sum_{l\mid m} 
\tfrac{1}{l}  \mu ( \tfrac{m}{l} ) 
\left(  \abs{\Hom(G,D_{2l})} -\abs{\Hom(G,\Z_{l})} \right)\\
 &=\sum_{l\mid m} 
\tfrac{1}{l}  \mu ( \tfrac{m}{l} ) 
\sum_{\rho \in \Epi(G,\Z_2)} \abs{Z^1_{\s\rho}(G,\Z_{l})} \\
\nonumber
& =\sum_{l\mid m} 
\tfrac{1}{l}  \mu ( \tfrac{m}{l} ) 
\sum_{\rho \in \Epi(G,\Z_2)}  \prod_{i=1}^{r} q_i^{d(\pi_{q_i}\s\rho)} .
\end{eqnarray}
In particular, $\abs{\Epi(F_n,D_{2m})} = (2^n-1) 
\sum _{ l \mid m} \frac{m}{l} \mu ( \frac{m}{l} ) l^n$, 
which, after some manipulations, recovers formula 
(\ref{eq:phi dihedral}).
} 
\end{example}

\section{Hall invariants and finite-index subgroups}
\label{sec:lis}
Let $G$ be a finitely generated group. For each positive integer $k$, let
$a_{k}(G)$ be the number of index $k$ subgroups of $G$. The 
behavior of the sequence $\{a_k(G)\}_{k\ge 1}$ (that is, the 
``subgroup growth" of $G$) has been the object of intense study 
ever since the foundational paper of M.~Hall~\cite{Ha}; see the 
monograph by A.~Lubotzky and D.~Segal \cite{LS} for a 
comprehensive survey.

Let $h_k(G)=\abs{\Hom(G,S_k)}$ and  $t_k(G)$ be the number 
of homomorphisms (respectively, transitive homomorphisms) 
from $G$ to the symmetric group $S_k$.   
It is readily seen that $a_{k}(G)=\frac{t_k(G)}{(k-1)!}$. 
The following recursion formula (due to M.~Hall \cite{Ha}) computes $a_k$ 
in terms of $h_1, \dots , h_k$, starting from  $a_1=h_1=1$:
\begin{equation}
\label{eq:mhall}
a_{k}(G)=\frac{1}{(k-1)!}\, h_{k}(G) -
\sum_{l=1}^{k-1} \frac{1}{(k-l)!}\, h_{k-l}(G) a_{l}(G).
\end{equation}

In this context, it is also useful to consider the 
{\em $\G$-Hall invariant} of $G$, 
\begin{equation}
\label{eq:delta hall}
\delta_{\Gamma}(G)=\abs{\Epi(G,\G)}/\abs{\Aut\Gamma}.
\end{equation}
Since $\Aut\Gamma$ acts freely and transitively on $\Epi(G,\Gamma)$,
the number $\delta_{\Gamma}(G)$ is an integer, which counts the homomorphic 
images of $G$ that are isomorphic to $\G$.  Notice that 
$\delta_{\Gamma_1\times \Gamma_2}(G)=\delta_{\Gamma_1}(G)
\delta_{\Gamma_2}(G)$, provided $\G_1$ and $\G_2$ have coprime orders. 

Formula (\ref{eq:mhall}), together with P.~Hall's 
enumeration principle (\ref{eq:hallenum1}), 
expresses the numbers $a_k=a_k(G)$ in terms of Hall invariants.  
For low indices, we have: $a_1=1$, $a_2=\delta_{\Z_2}$, 
$a_3=\delta_{\Z_3} + 3\delta_{S_3}$, and 
\begin{equation}
\label{eq:akdelta}
a_4=\tfrac{1}{2} \delta_{\Z_2}(1- \delta_{\Z_2}) 
+\delta_{\Z_4} + 4 \delta_{\Z_2^{\oplus 2}}  + 4 \delta_{D_8} 
+ 4 \delta_{A_4} + 4\delta_{S_4}.
\end{equation}

In general, the Hall invariants $\delta_{\G}(G)$ contain more 
information about a group $G$ than the numbers $a_k(G)$. 
For example, $a_k(\Pi_g)=a_k(\Pi_{2g}^*)$, for all 
$g\ge 1$ (see \cite{LS}), but clearly 
$\delta_{\Z_n}(\Pi_g)\ne \delta_{\Z_n}(\Pi_{2g}^*)$ 
for any odd $n>1$.  

When $G$ is finitely presented, all the Hall invariants  
that appear in (\ref{eq:akdelta}) can be expressed in 
terms of simple homological data.
If $\G$ is abelian, this is done in Theorem 3.1 in \cite{MS}. 
Let us briefly review how this goes.  

For a prime $p$, write the $p$-torsion part of $H_1(G,\Z)$ as 
$\bigoplus_{i\ge 1} \Z_{p^i}^{\oplus \alpha_i}$.  
Set $n=\rank H_1(G,\Z)$,  $\a=\sum i \alpha _i$, and 
$\b = \sum \alpha_i$.  
For a positive integer $s$, write 
$\a[s] = \sum_{i=1}^{s-1} i \alpha _i$. 
We then have:
\begin{align}
\label{eq:del abel}
\delta_{\Z_{p^{s}}}(G)
&=\tfrac{p^{sn+\a} -p^{(s-1)n+\a[s]}}{p^{s}-p^{s-1}},
\nonumber
\\
\delta_{\Z_{p}^{\oplus s}}(G)
&=\prod_{i=0}^{s-1}\tfrac{p^{n+\b}-p^{i}} {p^{s}-p^{i}},
\\
\delta_{\Z_{p}\oplus \Z_{p^{s}}}(G)&=
\tfrac{(p^{sn+\a} -p^{(s-1)n+\a[s]})(p^{n+\b}-p)}{p^{s+1}(p-1)^{2}}.
\nonumber
\end{align}
These formulas, together with the multiplicativity property 
of $\delta$ determine the $\Gamma$-invariants of $G$ for 
$\G$ abelian of order at most $31$, while higher orders are 
treated similarly. 

\renewcommand{\arraystretch}{2}
\begin{table}
\[
\begin{array}
{|c|c|}
\hline
\G &\delta_{\G}(G)
\\
\hline\hline
S_3 = \Z_3\rtimes_{\s} \Z_2 
&  \tfrac{1}{2} \sum_{\rho\in\Epi(G,\Z_2)} 
(3^{\dim_{\Z_3}H^{1}_{\s\rho}(G, \Z_3)}-1) 
\\
D_{8}=\Z_2\rtimes_{\chi} \Z_{2}^{\oplus 2} 
& \tfrac{1}{8} \sum_{\rho\in\Epi(G, \Z_2^{\oplus 2})} 
\epsilon_{\chi}(\rho) 2^{\dim_{\Z_2} H^1_{\rho}(G, \Z_2)} 
\\
Q_8 =\Z_2\times_{\chi} \Z_{2}^{\oplus 2}
& \tfrac{1}{24}\sum_{\rho\in\Epi(G, \Z_2^{\oplus 2})} 
\epsilon_{\chi}(\rho) 2^{\dim_{\Z_2} H^1_{\rho}(G, \Z_2)} 
\\
D_{12}=\Z_3\rtimes_{\s} \Z_{2}^{\oplus 2}
& \tfrac{1}{4} \sum_{\rho\in\Epi(G, \Z_2^{\oplus 2})}
\big( 3^{\dim_{\Z_3} H^1_{\s\rho}(G, \Z_3)} - 1 \big) 
\\
D^*_{12}=\Z_3\rtimes_{\s} \Z_4
& \tfrac{1}{4} \sum_{\rho\in\Epi(G, \Z_4)} 
(3^{\dim_{\Z_3} H^1_{\s\rho}(G, \Z_3)}-1)
\\
A_4=\Z_2^{\oplus 2} \rtimes_{\s} \Z_3 
&   \tfrac{1}{6} \sum_{\rho\in\Epi(G,\Z_3)} 
(2^{\dim_{\Z_2}H^{1}_{\s\rho}(G, \Z_2^{\oplus 2})}-1) 
\\
\hline
\end{array}
\]
\caption{\textsf{Hall invariants 
for non-abelian groups of order at most $12$.}}
\label{tab:hall12}
\end{table}
\renewcommand{\arraystretch}{1.1}

If $\G$ is non-abelian, of order at most $12$, the answer 
is given in Table \ref{tab:hall12}.  Plugging these answers, 
together with the ones from (\ref{eq:del s4}) and 
(\ref{eq:del abel} into formula (\ref{eq:akdelta})  
gives an expression for $a_4(G)$ solely in terms of 
cohomological invariants for $G$. 

Now let $a_{k}^{\nor}(G)$ be the number of index $k$, normal
subgroups of $G$.  Clearly,
\begin{equation}
\label{eq:aknor}
a_k^{\nor}(G)=\sum_{\abs{\Gamma}=k} \delta_{\Gamma} (G).
\end{equation}
In \cite{MS}, we used this formula to compute $a_k^{\nor}(G)$ in terms
of homological data, provided $k$ has at most two factors.  
Our approach worked for all $k\le 15$, except for $k=8$ and $k=12$.  
To compute $a_8^{\nor}$, we also needed to know $\delta_{D_8}$ 
and $\delta_{Q_8}$; for $a_{12}^{\nor}$, we also needed 
$\delta_{D_{12}}$ and $\delta_{D^*_{12}}$. The formulas in 
Table \ref{tab:hall12} complete the computation of 
$a_k^{\nor}(G)$, for $k\le 15$.

\section{Finite quotients of braid groups}
\label{sec:braids}

We conclude with a discussion of Artin's braid groups, viewed 
through the prism of their finite quotients and their finite-index subgroups.   
In addition to our own results, we use in crucial 
fashion results of Artin \cite{Ar},  Gorin and Lin \cite{GL}, and 
Lin \cite{Lin}, \cite{Lin2}, \cite{Lin3}. 

\subsection{Braid groups}
\label{subsec:bn}
The braid group on $n\ge 3$ strings has presentation
\[
B_n = \langle x, y \mid y ^n (y x)^{1-n}, \ 
[y ^ixy ^{-i},x], \ 2\le i\le n/2\rangle.
\]
Let $B_n'$ be the commutator subgroup. Clearly, 
$B_n/B_n'=\Z$, generated by $x$, and so we 
have a split extension, $B_n = B_n' \rtimes_{\tau} \Z$.  
It is also known that 
\begin{align*}
B_3& =F_2\rtimes_{\tau} \Z=\langle x, a,b \mid  a^x=b,\ 
b^x=b a^{-1} \rangle, \\
B_4 &=(F_2\rtimes F_2) \rtimes_{\tau} \Z= 
\left\langle \!
\begin{array}{l}
x, a,b \\ c, d 
\end{array}
\left|  
\begin{array}{llll}
a^x=b,
&b^x=b a^{-1},
&c^x= dc,
&d^x=d,
\\
c^a=d,
&c^b= d^{-1}c,
&d^a=dc^{-1}d^2,
&d^b= dc^{-1}d
\end{array}
\right. \!
\right\rangle \! .
\end{align*}
Note that $B'_3/B''_3=B'_4/B''_4=\Z^2$.  On the 
other hand, if $n\ge 5$, then $B'_n$ is perfect, see 
Gorin and Lin \cite{GL}. 

Now suppose $\G$ is a finite group. 
If $\G$ is cyclic,  then $\delta_{\G}(B_n)=1$.   
On the other hand, if $\G/\G'$ is non-cyclic, then $\delta_{\G}(B_n)=0$.  

If $n\ge 5$, and $\G$ is a finite quotient of $B_n$, then 
$\G'$ must be perfect.  Hence, every finite solvable quotient 
of $B_n$ must be cyclic, and so 
\begin{equation}
\label{eq:delbn}
\delta_{\G}(B_n) = \begin{cases} 
0 & \text{if $\G$ is not cyclic},\\
1 & \text{if $\G$ is cyclic},
\end{cases}
\end{equation}
whenever $\Gamma$ is a finite solvable group and  $n\ge 5$. 
On the other hand, the groups $B_3$ and $B_4$ have plenty 
of non-abelian, finite solvable quotients, as we see next.  

\subsection{Solvable quotients of $B_3$ and $B_4$}
\label{subsec:solvable b3b4}
Let $G$ be one of the braid groups $B_3$ or $B_4$.  
From the presentations above it is apparent that the 
maximal metabelian quotient, $G/G''$, is isomorphic to
$H=(\Z\oplus\Z)\rtimes_{\t} \Z$.  The monodromy action, 
$\tau=\big(\begin{smallmatrix} 0 & -1 \\ 1 & 1\end{smallmatrix}\big)$, 
has order $6$, and its characteristic polynomial is $t^2-t+1$. 

Now let $\G$ be a finite, metabelian quotient of $G$.  
Then $\G=\G'\rtimes_{\bar\t}\Z_k$, with $\G'$ a quotient of $\Z^2$. 
Write $\G'=\Z_m \oplus \Z_l$. 
We can pick generators $z\in \Z_k$ and 
$u, v \in \G'$  so that 
$\bar\t(z)\in \Aut(\G')$ is given by 
$\bar\t(u)=v$ and $\bar\t(v)=-u+v$. 
Hence, $u$ and $v$ have the same order, and 
so either $(m,l)=1$ or $l=m$.   
Analyzing the various possibilities, we obtain 
the following Proposition. 

\begin{proposition}
\label{prop:metabel b3b4}
Let $\G$ a finite, metabelian quotient of $G=B_3$ or $B_4$. 
Assume $\G$ is not cyclic. 
Then $\G$ is split metabelian, of  type
\begin{enumerate}
\item  
$\G=\Z_3\rtimes\Z_k$ with $k\equiv \pm 2 \md{6}$, 
in which case $\delta_{\G}(G)=1$; 
\item 
$\G=\Z_r\rtimes\Z_k$ with $r>3$ and $k\equiv 0 \md{6}$, 
in which case $\delta_{\G}(G)=2$; 

\item 
$\G=\Z_2^{\oplus 2}\rtimes\Z_k$ with $k\equiv 3 \md{6}$, 
in which case $\delta_{\G}(G)=1$; 

\item 
$\G=\Z_r^{\oplus 2}\rtimes\Z_k$ with $r>3$ and $k\equiv 0 \md{6}$, 
in which case $\delta_{\G}(G)=1$.  
\end{enumerate}
\end{proposition}

Now assume $\G$ is a finite, solvable, non-cyclic 
quotient of $B_3$ or $B_4$. If follows from the proof above that 
the maximal metabelian quotient, $\G/\G''$ has order 
divisible by $6$. But $\G$ is an extension of $\G/\G''$, 
and so $\abs{\G}$ is also divisible by $6$.  

Finite solvable quotients of $B_3$ can have 
derived length greater than $2$.  
For example, consider $S_4=\Z_2^{\oplus 2}\rtimes S_3$.  
We know $\abs{\Epi(B_3,S_3)}=6$. If $\rho$ is an 
epimorphism from $B_3$ to $S_3$, 
then $H^{1}_{\sigma\rho}(B_3;\Z_2^{\oplus 2})=\Z_2$. 
Hence, $\abs{\Epi(B_3,S_4)}=4\cdot 6 \cdot (2^{1} -1)$, 
and so $\delta_{S_4}(B_3) =1$. 

More generally,  if $\G_{r}=\Z_{2\cdot 3^r}\rtimes_{\chi} A_4$ 
is the sequence of groups starting from $\G_0=S_4$, 
then $\delta_{\G_r}(B_3)=1$. On the other hand, if 
$\tilde\G_r=\Z_{2\cdot 3^r}\rtimes_{\tilde\chi} A_4$ is the 
sequence of groups starting from $\tilde\G_0=\SL(2,3)$, then 
$\delta_{\tilde\G_r}(B_3)=2$.  

Since $B_4$ surjects onto $B_3$, it inherits all 
the finite quotients of $B_3$.  In general, though, $B_4$ has 
more epimorphisms onto a given finite 
quotient than $B_3$. The smallest solvable group  
for which this happens is $S_4$.  Indeed, 
$H^{1}_{\sigma\rho}(B_4;\Z_2^{\oplus 2})=\Z_2^{\oplus 2}$, 
for all  $\rho\colon B_4\surj S_3$; 
thus, $\abs{\Epi(B_4,S_4)}=4\cdot 6 \cdot (2^{2} -1)$, 
and so $\delta_{S_4}(B_4) =3$, although 
$\delta_{S_4}(B_3)=1$. 
 
In view of all this evidence, we propose the following conjecture. 
Let $\G$ be a finite solvable group.  Then 
\begin{equation}
\label{eq:delb3b4}
\delta_{\G}(B_3) \le 2 \quad \hbox{and}\quad 
\delta_{\G}(B_4) \le 3.
\end{equation}

\renewcommand{\arraystretch}{1.5}
{\small{
\begin{table}%[h]
\[
\begin{array}
{|c|c|c|c|c|c|c|c|c|c|c|c|c|c|c|c|}
\hline
   & a_3  & a_4 & a_5 & a_6 & 
a_7 & a_8 & a_9  & a_{10} & a_{11} & a_{12} 
& a_{13}  & a_{14} & a_{15} & a_{16} 
\\ 
\hline
B_3 & \mathbf{4} &  \mathbf{9} & \mathbf{6} & \mathbf{22} & \mathbf{43} 
& 49& 130& 266& 287& 786 & 1730 & 2199 & 5184 &12193 
\\
\hline
B_4 &  \mathbf{4} &  \mathbf{17} & \mathbf{6} & \mathbf{34} 
& 43& 81& 148& 266& 287& 938& 1730
& 2199& 5199& 12449 
\\
\hline
B_5 &  \mathbf{1}&  \mathbf{1}& \mathbf{6} & \mathbf{7}&  \mathbf{1} 
& 1& 1& 26& 1& 19 &  1 & 1 & 36 & 17 
\\
\hline
B_6  &  \mathbf{1}&  \mathbf{1}&  \mathbf{1}& \mathbf{13}&  \mathbf{1}&  
\mathbf{1}&  \mathbf{1}& \mathbf{11} & 1& 25 &1 & 1 & 31 & 1  
\\
\hline
B_7 &  \mathbf{1}&  \mathbf{1}&  \mathbf{1}&  \mathbf{1}&  \mathbf{8}&  
\mathbf{1}&  \mathbf{1}&  \mathbf{1}&  \mathbf{1}&  \mathbf{1}&  
\mathbf{1}& \mathbf{22} & 1 & 1
\\
\hline
B_8 &  \mathbf{1}&  \mathbf{1}&  \mathbf{1}&  \mathbf{1}&  \mathbf{1}&  
\mathbf{9}&  \mathbf{1}&  \mathbf{1}&  \mathbf{1}&  \mathbf{1}&  
\mathbf{1} &  \mathbf{1} &  \mathbf{1} & \mathbf{25}
\\
\hline
\end{array}
\]
\caption{\textsf{Number of low-index subgroups of $B_n$ ($n\le 8$)}}
\label{tab:akbn}
\end{table}
}}

\subsection{Finite-index subgroups of $B_n$}
\label{subsec:akbn}
We conclude with a discussion of the subgroup growth 
of the braid groups $B_n$.  Of course, $a_1(B_n)=a_2(B_n)=1$. 
The values of $a_k(B_n)$ for $3\le n\le 8$ and $3\le k\le 16$ 
are listed in Table~\ref{tab:akbn}.  The values not in bold 
were computed solely by machine.  The values in bold can  
be justified, as follows.  

\begin{enumerate}
\item \label{b1}
Using results from \S \ref{subsec:solvable b3b4}, 
we see that $a_3(B_3)=a_3(B_4)=4$, 
and $a_3(B_n)=1$, for $n>4$. Furthermore, 
$a_4(B_3)=9$, $a_4(B_4)=17$, and $a_4(B_n)=1$, for $n>4$.

\item \label{b2}
Proposition 4.1 from Lin \cite{Lin3} give 
$t_k(B_3)$, for $4\le k\le 7$, while 
Propositions 4.4 and 4.7 from \cite{Lin3}  
give $t_5(B_4)$ and $t_6(B_4)$. 
This gives the corresponding values 
for $a_k(B_3)$ and $a_k(B_4)$.

\item \label{b3}
Suppose $n>4$ and $k<n$. In \cite{Lin}, Lin showed that 
any transitive homomorphism $B_n \to S_k$ 
has cyclic image. This implies $t_k(B_n)=(k-1)!$, and so 
$a_k(B_n)=1$. 

\item \label{b4}
In \cite{Ar}, Artin computed $\abs{\Epi(B_n, S_n)}$ for all $n$. 
This gives $a_5(B_5)=6$, $a_6(B_6)=13$, and 
$a_n(B_n)=n+1$, for $n>6$.  

\item \label{b5}
Suppose $6<n<k<2n$. In Theorem F.a) from \cite{Lin3}, Lin proves 
that any transitive homomorphism  $B_n \to S_k$  has cyclic image.
Consequently, $t_k(B_n)=(k-1)!$ and so $a_k(B_n)=1$. 

\item \label{b6}
Suppose $n>6$. Up to conjugation, there are $4$ 
transitive homomorphisms $B_n\to S_{2n}$, of which 
$3$ are non-cyclic, see \cite[Theorem F.b)]{Lin3}. 
It is readily seen that the centralizer of those $3$ homomorphisms 
is the involution $(1,2)(3,4)\cdots(2n-1,2n)$.  Hence, 
$t_{2n}(B_n)=(2n-1)! +3(2n)!/2$, and so $a_{2n}(B_n)=3n+1$. 

\item \label{b7}  
Further results from Sections 4 and 7 in \cite{Lin3} give 
$t_6(B_5)$, $t_7(B_5)$, and $t_k(B_6)$, for $7\le k\le 10$;   
the corresponding values for $a_k(B_n)$ follow. 
\end{enumerate}

%\begin{itemize}
%\item Proposition 4.10 gives $t_6(B_5)$;
%\item Theorem Eb) and Proposition 7.17 give $t_7(B_5)$;
%\item Theorem Ea) and Theorem 7.4 give $t_7(B_6)$;
%\item Theorem 7.20 and Proposition 7.19 give $t_8(B_6)$;
%\item Proposition 7.23 gives $t_9(B_6)$;
%\item Remark 7.24 gives $t_{10}(B_6)$;
%\end{itemize}

Out of this discussion, we obtain the following corollary.
\begin{corollary}
\label{cor:akbn}
For the specified values of $k$ and $n$, the number of 
index $k$ subgroups of the braid group $B_n$ is given by
\[
a_k(B_n)=
\begin{cases}
1 & \text{for $k<n$ and $n>4$, or $6< n< k < 2n$}, \\
n+1 & \text{if $k=n$ and $n>6$}, \\
3n+1 & \text{if $k=2n$ and $n>6$}.
\end{cases}
\]
\end{corollary}

Based on this evidence, we propose the following conjecture.

\begin{conj}
\label{conj:akbn}
For all $n\gg 0$, 
\[
a_k(B_n)= 
\begin{cases}
1 & \text{if $n\nmid k$}, \\
c(k/n)\cdot n+1 & \text{if $n\mid k$}
\end{cases}
\]
where $c(k/n)$ is a constant, depending only on $k/n$. 
\end{conj}

\medskip
\begin{ack}
A substantial amount of this work was done while the first author 
was visiting Northeastern University in February-May, 2004, with 
support from the Research and Scholarship Development Fund 
and the Mathematics Department.  

The authors are grateful to Vladimir Lin for sending a preliminary
version of \cite{Lin3}, and for his comments on the results therein. 
They also thank the referee for pointing out an inaccuracy in the 
original statement of Theorem \ref{thm:gg}, and for many helpful 
suggestions.  

Thanks to Tony Iarrobino for providing a brand new 
PowerMac G4 for this project. The computations were 
carried out with the help of the package {\sl GAP~4.4} \cite{gap}.  

\end{ack}

\end{document}